\documentclass[a4paper]{amsart}
\usepackage[utf8]{inputenc}
\usepackage{eucal}
\usepackage{amssymb,amsmath,amsthm}
\usepackage{lmodern}
\usepackage[T1]{fontenc}
\usepackage[hidelinks]{hyperref}
\usepackage{mathrsfs}
\usepackage{caption}
\usepackage{mathtools}

\allowdisplaybreaks

\usepackage{thmtools}

\declaretheorem[name=Remark, style=definition, qed=$\triangleleft$, numbered=no]{rem}

\newcommand{\mc}{\mathcal}
\newcommand{\on}{\operatorname}

\newcommand{\la}{\langle}
\newcommand{\ra}{\rangle}
\newcommand{\R}{\mathbb{R}}

\title{Ricci tensor in graded geometry}

\author{Fridrich Valach}
\address{Mathematical Institute, Faculty of Mathematics and Physics, Charles University Prague 186 75, Czech Republic}
\email{fridrich.valach@gmail.com}
\thanks{This work was supported by the GA\v{C}R Grant EXPRO 19-28628X}

\begin{document}

\begin{abstract}
We define the notion of the Ricci tensor for NQ symplectic manifolds of degree 2 and show that it corresponds to the standard generalized Ricci tensor on Courant algebroids. We use an appropriate notion of connections compatible with the generalized metric on the graded manifold.
\end{abstract}

\maketitle

\section{Introduction}
It has been known for some time that Courant algebroids \cite{lxw} provide a natural framework for the study of several aspects of string theory. First, the string sigma models can be seen \cite{SCS} as particular cases of the Courant sigma models \cite{I,roy} on manifolds with boundary. Furthermore, the low-energy dynamics of string theory can be described via a suitable generalization of the Ricci tensor and scalar curvature, defined for Courant algebroids \cite{MG,CSCW,MGF2,MGF,JV,SVsugra}. 

While for some purposes it is enough to consider only a special class of the so-called \emph{exact} Courant algebroids \cite{let}, for the study of dualities (such as the Poisson-Lie T-duality \cite{KS}), it is neccessary to understand the constructions in the general setup \cite{let}. Such generalized Ricci tensor was introduced and studied in \cite{MGF,JV,SVsugra}. In particular, its properties were used to provide a proof of the compatibility of the Poisson-Lie T-duality with the RG flow and with the string background equations, extending the results of \cite{VKS} and \cite{H}, as well as to find new solutions to generalized supergravity equations.

On the other hand, it is known \cite{let,DR} that graded geometry offers a much more conceptual viewpoint of Courant algebroids, simplifying the formulas and providing various new insights. It is therefore desirable to also find a formulation of the curvature tensors in the graded language.

Motivated by the recent work \cite{A}, where the graded analogues of connections, curvature and torsion were introduced and studied, we propose a simple definition of the generalized Ricci tensor in the graded setup. The purpose of this is twofold -- it provides a more conceptual viewpoint of the generalized Ricci tensor and at the same time opens a very concrete door towards generalizations, for example in the context of U-duality, following \cite{CSCW2}.  This will be explored in a future work.

This note is structured as follows. We start by reviewing the notion of generalized metric, from the perspective of graded geometry. We then recall the neccessary definitions from \cite{A}, introduce the notion of connections with invariant torsion and present the Ricci tensor. Finally, we explore the relation to the more standard Courant algebroid connections \cite{MG} and discuss the exact case, where this Ricci tensor produces the usual one from Riemannian geometry.

\subsection*{Acknowledgements}
The author would like to thank Andreas Deser and Pavol \v Severa for helpful discussions, suggestions, and comments on the preliminary version of the paper. The author would also like to acknowledge the COST action MP 1405 Quantum structure of spacetime and the Corfu Summer Institute 2019 at EISA, where the main idea of the paper was conceived.

\section{Generalized metric}
In what follows, we shall use $\mc E$ to denote an NQ symplectic manifold of degree 2 \cite{some}. This means that $\mc E$ is an $\mathbb{N}$-graded manifold\footnote{i.e.\ an ordinary manifold $\mc E_0$ together with a sheaf of $\mathbb{N}$-graded commutative algebras, locally of the form $C^\infty(\mc U)\otimes S(V)$, with $S(V)$ the free graded commutative algebra generated by a finite-dimensional vector space $V=\bigoplus_{i=1}^mV_m$, and $\mc U$ an open subset of $\mc E_0$}, equipped with a symplectic form of degree 2, and a degree 1 symplectic vector field $Q_\mc E$, satisfying $Q^2_\mc E=0$. There is an associated sequence of fibrations
\[\mc E=\mc E_2 \to \mc E_1 \to \mc E_0,\]
which corresponds to the subsheafs generated by coordinates of degrees up to 2, up to 1, and up to 0, respectively. In particular, the last arrow gives a vector bundle.

A \emph{generalized metric} is, in the graded language, a symplectic involution $\iota$ on $\mc E$, which preserves the basis $\mc E_0$, i.e. it is a diffeomorphism of $\mc E$ satisfying
\[\iota^2=id_\mc E,\quad \iota^*\omega=\omega,\quad \iota|_{\mc E_0}=id_{\mc E_0}.\]

Given a generalized metric, it is always possible to locally choose the coordinates $x^i$, $e^a$, $e^{\dot a}$, $p_i$ on $\mc E$ of degrees 0, 1, 1, 2, respectively, such that
\[\omega=dp_idx^i+de_a de^a +de_{\dot a}de^{\dot a},\]
\[\iota^*x^i=x^i,\quad \iota^*p_i=p_i,\quad \iota^*e^a=e^a,\quad \iota^*e^{\dot a}=-e^{\dot a},\]
where $e_a:=g_{ab}e^b$, $e_{\dot a}:=g_{\dot a\dot b}e^{\dot b}$ for some constants $g_{ab}$, $g_{\dot a\dot b}$. We will sometimes denote the coordinates $e^a$, $e^{\dot a}$ collectively as $e^\alpha$.\footnote{The proof goes as follows: Using the (graded) Darboux theorem we first find coordinates $x^i$, $e^\alpha$, $p_i$ such that $\omega=dp_idx^i+g_{\alpha\beta}de^\alpha de^\beta$, for $g_{\alpha\beta}$ a diagonal matrix with only 1 and $-1$ on the diagonal. It follows that the ($x$-dependent) matrix $R^\alpha_\beta$ defined by $\iota^*e^\alpha=R^\alpha_\beta(x) e^\beta$ is idempotent and orthogonal w.r.t.\ $g$, and thus can be made into the form $\on{diag}(1,\dots,1,-1,\dots,-1)$ using $e'^\alpha=O^\alpha_\beta(x) e^\beta$, for $O$ orthogonal. Making an appropriate shift of $p_i$, the form of $\omega$ is preserved in the new coordinates.} 

It is easy to see that $Q_\mc E$ is always Hamiltonian, i.e.\ it comes from a degree 3 function $H$. The most general such function has the form
\[H=\rho_\alpha^i(x)p_ie^\alpha-\tfrac16 c_{\alpha\beta\gamma}(x)e^\alpha e^\beta e^\gamma,\vspace{-.2cm}\]
giving 
\[Q_\mc E=\rho^i_\alpha e^\alpha\partial_{x^i}+(\rho^i_\alpha p_i-\tfrac12 c_{\alpha\beta\gamma}e^\beta e^\gamma)\partial_{e_\alpha}+(\tfrac16 c_{\alpha \beta\gamma,i}e^\alpha e^\beta e^\gamma-\rho^j_{\alpha,i}e^\alpha p_j)\partial_{p_i},\]
where, in setting $e_\alpha:=g_{\alpha\beta}e^\beta$, we extended $g_{ab}$, $g_{\dot a\dot b}$ to $g_{\alpha\beta}$ by adding $g_{a\dot b}=g_{\dot a b}=0$.
The condition $Q^2_\mc E=0$ translates to the classical master equation $\{H,H\}=0$.

Since $\iota$ preserves the degree, it induces an involution on $\mc E_1$. Because of the base-fixing condition, $\iota$ is in turn fully determined by the corresponding fixed point set $\mc V_1\subset \mc E_1$. Seeing $\mc V_1$ as a vector bundle over $\mc E_0$, we can pull it back along $\mc E\to\mc E_0$ to obtain a bundle $\mc V\to\mc E$. The latter bundle is locally described by coordinates $x^i$, $e^a$, $\xi^a$, $e^{\dot a}$, $p_i$, with one copy of $\xi^a$, $\deg\xi^a=1$, corresponding to each $e^a$.
Finally, notice that $\mc E_1$, $\mc V_1$, and $\mc V$ are all symplectic vector bundles.

\begin{rem}
The relation to Courant algebroids is as follows \cite{DR,let}. The Courant algebroid is given by the ordinary vector bundle $E=\mc E_1[-1]$ over $\mc E_0$. (In other words, $x^i$ are coordinates on the base and $e^\alpha$ correspond to the linear coordinates on the fibers of the algebroid.) The coefficients $c$ and $\rho$ give the structure functions of the bracket and the anchor, respectively, while $g$ encodes the fiberwise inner product. Finally, the involution $\iota$ corresponds to the usual viewpoint of generalized metric as a fiberwise reflection on the Courant algebroid, or equivalently, as a subbundle $V_+=\mc V_1[-1]\subset E$.
\end{rem}

\section{Tautological section and contraction}

Let us now denote the vector bundle morphism $\mc V\to \mc V_1$ by $\varphi$. There is a unique section $\tau'\colon \mc E\to\mc V$ such that $\varphi\circ \tau'$ coincides with the map $\mc E\to \mc E_1 \to \mc V_1$ (the last arrow is the orthogonal projection). Since the bundle $\mc V$ is symplectic, we get an induced section $\tau$ on the dual bundle $\mc V^*\to \mc E$. We call $\tau$ the \emph{tautological section} \cite{A}. More concretely, identifying sections of $\mc V^*$ with functions on $\mc V$ which are linear in the fiber coordinates, we get
\[\tau=e_a\xi^a.\]

The involution on $C^\infty(\mc V)$ (induced by $\iota$) allows us to split any vector field on $\mc V$ into the sum of its self-dual and anti-self-dual part. We will denote the anti-self-dual part of a vector field $D$ by $\pi D$. Let us now consider a special subspace $\text{End}_2$ of the space of vector fields, given by degree 2 bundle endomorphisms of $\mc V^*$. Locally we have
\[D=D^{a\phantom{b}i}_{\phantom{a}b}(x)p_i\xi^b\partial_{\xi^a}+\tfrac12 D^a_{\phantom a b \alpha\beta}(x)e^\alpha e^\beta\xi^b\partial_{\xi^a}\quad \xmapsto{\pi}\quad D^a_{\phantom a b c\dot{d}}(x)e^c e^{\dot{d}}\xi^b\partial_{\xi^a}.\]
Writing $\mc V_1^\perp\subset \mc E_1$ for the subbundle (over $\mc E_0$) perpendicular to $\mc V_1\subset\mc E_1$, we have an identification (we use $\Gamma$ for the space of sections) \[\pi(\text{End}_2)\cong \Gamma(\mc V_1^*\otimes \mc V_1^{\perp*}\otimes End(\mc V_1^*))\cong \Gamma(\mc V_1^*\otimes \mc V_1^{\perp*}\otimes \mc V_1\otimes \mc V_1^*).\]
We define the \emph{contraction map} $C\colon \text{End}_2\to \Gamma(\mc V_1^{\perp *}\otimes \mc V_1^*)$ as the projection $\pi$ followed by the contraction of the first and third factor in the last expression.
Explicitly,
\[C\colon D\mapsto \partial_{e^a} (\pi D) \xi^a.\]

\section{Connections, torsion and curvature}
Following \cite{A}, a \emph{connection} on a $\mc V^* \to \mc E$ is a degree 1 vector field $Q$ on $\mc V$, which projects to $Q_\mc E$, and which preserves the space of sections $\Gamma(\mc V^*)\subset C^\infty(\mc V)$. 
The \emph{torsion} is then a particular section of $\mc V^*$, defined as $Q\tau$. The \emph{curvature} of $Q$ is the vector field $Q^2\equiv \tfrac12[Q,Q]$ on $\mc V$. One easily sees that $Q^2\in\text{End}_2$.

We will say that a connection $Q$ on $\mc V$ \emph{has invariant torsion} if its torsion is invariant under the induced involution on $\Gamma(\mc V^*)$.
Finally, we define the \emph{Ricci tensor} $\on{Ric}$ by
\begin{equation}\label{ric_def}
\on{Ric}:=CQ^2\in \Gamma(\mc V_1^{\perp *}\otimes \mc V_1^*).
\end{equation}

\begin{rem}
Connections with invariant torsion should be seen as analogues of the usual Levi-Civita connection -- the compatibility with the metric is replaced by the fact that $Q$ is a vector field on $\mc V$ and the vanishing of the torsion is replaced by the invariance of $Q\tau$ (as seen below, in general it is impossible to have $Q\tau=0$; on the other hand the invariant torsion condition does not fix the connection uniquely, c.f.\ \cite{CSCW}).

The contraction in (\ref{ric_def}) corresponds to the usual procedure for obtaining the Ricci tensor from the Riemann tensor. The insertion of the projection $\pi$ keeps only the part of the tensor which can be identified with an infinitesimal deformation of the generalized metric.\footnote{It also allows us to define the contraction.} This is due to the identification 
\[\Gamma(\mc V_1^{\perp *}\otimes\mc V_1^*)\cong\on{Hom}_{\mc E_0}(\mc V_1,\mc V_1^{\perp *})\cong\on{Hom}_{\mc E_0}(\mc V_1,\mc V_1^{\perp}),\] with the last isomorphism provided by the symplectic form on $\mc V_1^{\perp}$. This deformation is the (infinitesimal) \emph{generalized Ricci flow} \cite{SV,MGF}.
\end{rem}

Passing again to a coordinate description, a general connection on $\mc V^*$ has the form
\[Q=Q_\mc E+\psi^a_{\phantom{a}b\alpha }(x)e^\alpha \xi^b\partial_{\xi^a}.\]
Its torsion is 
\[Q\tau=\rho_a^i p_i\xi^a-\tfrac12 c_{a\beta\gamma}e^\beta e^\gamma\xi^a-\psi_{\phantom{a}b\alpha}^a e^\alpha\xi^b e_a.\]
Since the involution preserves the coordinates $\xi^a$, the invariance of torsion is equivalent to the constraint
\begin{equation}
\label{eq:lc}
\psi_{\phantom{a}b\dot a}^a=c_{\phantom{a}b\dot a}^a.
\end{equation}
In particular, the coefficients $\psi^a_{\phantom{a}bc}$ are left unrestricted. However, as we will see, in the case of invariant torsion the curvature only depends on $\psi^a_{\phantom{a}bc}$ through the trace \[\lambda_b:=\psi^a_{\phantom{a}ba}.\]

Concretely, a short calculation (see Appendix) reveals that for a connection $Q$ with invariant torsion,
\begin{equation}\label{ric}
\on{Ric}=(c^c_{\phantom{c}b\dot{a} }\lambda_c-\rho^i_a c^a_{\phantom{a}b\dot{a},i }
  +\rho^i_{\dot{a}} \lambda_{b,i}
  +c_{\dot{c}a\dot{a}}c^{a\phantom{b}\dot{c}}_{\phantom{a}b}
  )\xi^b e^{\dot a} .
\end{equation}

This recovers exactly the formula for the generalized Ricci tensor from \cite{SVsugra}.

\begin{rem}
More precisely, using the notation of \cite{SVsugra} we have 
\[\on{Ric}=\on{GRic}_{V_+,\on{div}}(e_c,e_{\dot{a}})\xi^c e^{\dot{a}},\]
where the divergence operator is given by $\on{div}(e_a)=-\lambda_a$.
\end{rem}

\section{Connection with Courant algebroid connections}
From the definition it follows that connections are in one-to-one correspondence with linear maps\footnote{We use subscript (for $C^\infty$ and $\Gamma$) to denote the degree.} $\hat Q\colon \Gamma_1(\mc V^*) \to \Gamma_2(\mc V^*)$, satisfying $\hat Q(f u)=f(\hat Qu)+(Q_\mc E f)u$, for $f\in C^\infty_0(\mc E)$, $u\in\Gamma_1(\mc V^*)$. Since
\[C^\infty_0(\mc E)\cong C^\infty(\mc E_0),\quad C^\infty_1(\mc E)\cong \Gamma(E^*),\quad \Gamma_1(\mc V^*)\cong \Gamma(V_+^*),\quad \Gamma_2(\mc V^*)\cong \Gamma(E^*\otimes V_+^*).\]
we can understand $\hat Q$ as a map
\[\nabla\colon \Gamma(V_+^*)\to\Gamma(E^*\otimes V_+^*),\quad\text{such that}\quad \nabla(fu)=f\nabla u+(Q_\mc E f)\otimes u.\]
Dually, we have a map $\nabla\colon \Gamma(V_+)\to\Gamma(E^*\otimes V_+)$ (satisfying the same Leibniz identity). This is known in the literature as the \emph{Courant algebroid (or generalized) connection} \cite{MG}.

If $Q$ has invariant torsion, then it is uniquely determined (c.f.\ (\ref{eq:lc})) by the restriction of $\nabla$ to a map
$\Gamma(V_+)\to\Gamma(V_+^*\otimes V_+)$.

\section{The exact case}
Let us now consider the following example \cite{let,SCS}.\footnote{Every $\mc E$ with vanishing tangential cohomology and a generic generalized metric can be put into this form \cite{let,MG2}.} First,
\[\mc E=T^*[2]T[1]\mc E_0, \quad H=d+\eta, \quad \text{ for }\eta\in\Omega^3_{\text{closed}}(\mc E_0),\]
with the standard symplectic form on the cotangent bundle. Here $d$ is understood as a vector field on $T[1]\mc E_0$ (thus a function on $T^*[2]T[1]\mc E_0$), and $\eta\in \Omega(M)\cong C^\infty(T[1]M)$ is pulled back to $T^*[2]T[1]\mc E_0$. (Such $\mc E$ is called \emph{exact}.) We fix the generalized metric by the requirement that the submanifold $\mc V_1\subset \mc E_1\cong (T\oplus T^*)[1]\mc E_0$ corresponds to the graph of a (pseudo Riemannian) metric $\la\cdot,\cdot\ra$ on $\mc E_0$. 

In particular, we have a vector-bundle isomorphism $T\mc E_0\cong V_+$. Combining this with the result of the previous section, a connection on $\mc V^*$ with invariant torsion is now given by an ordinary\footnote{The term ``ordinary'' refers here to the usual non-graded geometry.} connection on $T\mc E_0$. Let us now take $Q$ given by the ordinary Levi-Civita connection on $T\mc E_0$ w.r.t.\ the metric on $\mc E_0$.

Choose a local frame $E_a$ on $T\mc E_0$ satisfying $\la E_a,E_b\ra=\pm \delta_{ab}$. We define the frame $F_a$ on $T^*\mc E_0$ by $F_a:=\la E_a,\cdot\ra$ and we denote by $E^a$, $F^a$ the induced fiber coordinates on $\mc E_1\cong (T\oplus T^*)[1]\mc E_0$. Finally, setting $e^a:=\tfrac12 (E^a+F^a)$, $e^{\dot b}:=\tfrac12 (E^b-F^b)$, we have
\[\on{Ric}=\on{Ric}_{\eta}(E_a,E_b)e^a e^{\dot b},\]
where $\on{Ric}_\eta$ is the ordinary Ricci tensor on $T\mc E_0$ for the metric connection (w.r.t.\ $\la\cdot,\cdot\ra$) with torsion given by $\eta$. For the proof of this fact we refer the reader to \cite{SVsugra}. It can also be verified by a direct calculation.

\appendix
\section{Expression for the Ricci tensor}
Using the fact that $Q^2_\mc E=0$, we have (for a general $Q$)
\begin{align*}
  Q^2&=(Q_\mc E \psi^a_{\phantom{a}b\alpha }e^\alpha \xi^b - \psi^b_{\phantom{b}c\beta }e^\beta \xi^c\psi^a_{\phantom{a}b\alpha }e^\alpha)\partial_{\xi^a}\\
  &=(\rho^i_\beta e^\beta \psi^a_{\phantom{a}b\alpha,i }e^\alpha \xi^b
  +\rho^i_\alpha p_i \psi^{a\phantom{b}\alpha}_{\phantom{a}b} \xi^b
  -\tfrac12 c_{\alpha\beta\gamma}e^\beta e^\gamma \psi^{a\phantom{b}\alpha}_{\phantom{a}b}\xi^b
  -\psi^b_{\phantom{b}c\beta }e^\beta \xi^c\psi^a_{\phantom{a}b\alpha }e^\alpha)\partial_{\xi^a}\\
\pi Q^2 &= (\rho^i_d \psi^a_{\phantom{a}b\dot{a},i}e^de^{\dot{a}}\xi^b+\rho^i_{\dot{a}}\psi^a_{\phantom{a}bd,i}e^{\dot{a}}e^d\xi^b-c_{\alpha d \dot{a}}\psi^{a\phantom{b}\alpha}_{\phantom{a}b}e^de^{\dot{a}}\xi^b+\psi^b_{\phantom{b}cd}\psi^a_{\phantom{a}b\dot{a}}e^de^{\dot{a}}\xi^c\\
&+\psi^b_{\phantom{b}c\dot{a}}\psi^a_{\phantom{a}bd}e^{\dot{a}}e^d\xi^c)\partial_{\xi^a}\\
CQ^2&=(\rho^i_a\psi^a_{\phantom{a}b\dot{a},i}-\rho^i_{\dot{a}}\psi^a_{\phantom{a}ba,i}-c_{\alpha a\dot{a}}\psi^{a\phantom{b}\alpha}_{\phantom{a}b}+\psi^c_{\phantom{c}ba}\psi^a_{\phantom{a}c\dot{a}}-\psi^c_{\phantom{c}b\dot{a}}\psi^a_{\phantom{a}ca})e^{\dot{a}}\xi^b
\end{align*}
Using $\psi^a_{\phantom{a}b\dot{a}}=c^a_{\phantom{a}b\dot{a}}$ and $\psi^a_{\phantom{a}ba}=\lambda_b$ we recover (\ref{ric}).

\end{document}